\documentclass{amsart}
\usepackage{graphicx}
\usepackage{amssymb,latexsym,amsmath,mathrsfs} 
\usepackage{epsfig}
\usepackage{color}

\newcommand{\RR}{\mathbb{R}}

\def\BibTeX{{\rm B\kern-.05em{\sc i\kern-.025em b}\kern-.08em
    T\kern-.1667em\lower.7ex\hbox{E}\kern-.125emX}}

\newtheorem{thm}{Theorem}[section]

\newtheorem{exmp}[thm]{Example}

\theoremstyle{definition}

\theoremstyle{remark}

\numberwithin{equation}{section}

\begin{document}

\title []{A line search algorithm for Wind field adjustment with incomplete data and RBF approximation}  

%
\author[]{Daniel A. Cervantes}
\address[Daniel A. Cervantes]{Instituto de Matem\'aticas \\
        UNAM\\
        CDMX.}
\email[D.~Cervantes]{dcchivela@ciencias.unam.mx}

\author[]{Pedro Gonz\'alez Casanova}
\address[Pedro Gonz\'alez Casanova]{Instituto de Matem\'aticas \\
        UNAM\\
        CDMX.}
\email[P.~Gonzalez]{casanova@matem.unam.mx}

\author[]{Christian Gout}
\address[Christian Gout]{INSA Rouen \\
        France.}
        
\email[Christian Gout]{christian.gout@insa-rouen.fr}

\author[]{Miguel \'Angel Moreles*}
\address[Miguel \'Angel Moreles]{Centro de Investigaciones en Matem\'aticas\\
        Guanajuato M\'exico}
 
\email[M. A.~Moreles]{moreles@cimat.mx}

\thanks{*Corresponding author}

\date{}


\maketitle

\begin{abstract}
The problem of concern in this work is the construction of free divergence fields given scattered horizontal components. As customary, the problem is formulated as a PDE constrained least squares problem. The novelty of our approach is to construct the so called adjusted field, as the unique solution along an appropriately chosen descent direction. The latter is obtained by the adjoint equation technique. It is shown that the classical adjusted field of Sasaki's is a particular case.  On choosing descent directions, the underlying mass consistent model leads to the solution of an elliptic problem which is solved by means of a Radial Basis Functions method. Finally some numerical results for wind field adjustment are presented.
\end{abstract}

\section{Introduction}
\label{intro}
The problem of recovering atmospheric wind fields from prescribed horizontal data is of great interest in meteorological applications.  In practice the vertical component is unavailable. Consequently, measured data is complemented with mass consistency to pose a variational problem for approximation. This approach dates back to Sasaki \cite{sasaki}. Literature on the subject is vast, a timely review is presented in Ratto et al \cite{ratto} .

The numerical approximation of the variational problem requires the solution of Poisson boundary value problems. Numerical approximations of the solutions can be obtained conventionally using the Finite Element Method, Finite Volume Method, Finite Differences, etc. In these methods mesh managing is computationally expensive. For the wind field adjustment problem, a case can be made for a mesh-free approach in terms of Radial Basis Functions (RBF). See Pepper, Rasmussen \& Fyda \cite{Pepperetal}, and Cervantes et al \cite{cervantesetal}.

In the context of RBF,  the problem of wind field recovery has been also considered as a smoothing problem. Mass consistency is introduced by penalizing with the norm of the divergence of the vector field. For the case of polyharmonic splines, Benbourhim and Bouhamidi \cite{najib}, developed a smoothing algorithm for a given set of prescribed data in two and three dimensions. A full convergence analysis is provided. These RBF techniques can  be traced back to Duchon's works on thin plate splines \cite{duchon}. 

The objective here, is to build on these methodologies and contribute on two unresolved issues. The wind field adjustment problem is posed on a bounded domain and boundary conditions are to be imposed, for instance the topography of the surface. The smoothing approach considers for approximation a linear combination of polyharmonic splines, imposing boundary conditions is not straightforward.  In the Sasaki's approach, boundary conditions are imposed somewhat heuristically. Also, given that only the first two components of the field are known, an initial guess of the third component is required in this approach. Common practice is unsatisfactory: \emph{ since measurements of the vertical velocity component are seldom available, the initial vertical is usually set to zero.}  Ratto et al \cite{ratto}.

Consequently, in this work the problem is dealt with the data available in practice. Namely,  the least squares functional to be introduced only involves the two known components. Minimization is carried out on descent directions, boundary conditions arise naturally and can be imposed on physical grounds. A RBF method is used for approximation.  Finally, a proper abstract setting for analysis and numerics is provided.

This article is organized as follows.

In Section 2 the Hilbert function spaces are introduced, and the variational problem is formulated. From basic results form minimization of convex functionals on Hilbert spaces, existence of minima is proven. The construction of the adjusted mass consistent field is the content of Section 3. Also, the Sasaki solution is derived as a particular case. In Section 4 we discuss boundary conditions for the Poisson problems arising form the line search algorithm, and the RBF discretisation. Numerical results are presented in Section 5, noteworthy the vortex field in  Benbourhim and Bouhamidi. Conclusions and comments on future research close this exposition.

\section{Reconstruction of mass consistent vector fields with unknown vertical component}

\subsection{Problem Formulation}

Assume the first two components of a $3D$ vector field are known at $N$ non uniform points. Namely $\mathbf{U}^0_i=(u_1(\mathbf{x}_i),u_2(\mathbf{x}_i)),\ i=1,2,\ldots,N$, for $\mathbf{x}_i\in\Omega\subset\mathbb{R}^3$.

\bigskip

\noindent\textbf{Problem.} Construct a flow field $\mathbf{u}(\mathbf{x})=(u_1(\mathbf{x}),u_2(\mathbf{x}),u_3(\mathbf{x}))$  in the entire domain $\Omega$, assuming that the set of discrete field values are known and that the approximant satisfies the continuity equation
\[
\nabla\cdot\mathbf{u}=0.
\]

The classical approach is as follows:

\begin{enumerate}
\item Initialization. From $\mathbf{U}^0_i,\ i=1,2,\ldots,N$, determine an initial field $\mathbf{u}^0(x)$ entirely in $\Omega$.

\item Correction. By means of a minimization procedure, adjust $\mathbf{u}^0$ to a reconstructed field satisfying mass conservation.

\end{enumerate}

\smallskip

Our focus is on (2).

\bigskip

\noindent\textbf{Function spaces}

Let $\Omega\subset\mathbb{R}^ 3$ with boundary $\Gamma$. We assume that $\Omega$ is a Lipschitz domain. We shall be concerned with $n-$dimensional vector functions. As customary we denote
\[
\mathbf{X}=X^n,
\]
e.g.,
\[
\mathbf{L}^2(\Omega)=\left( L^2(\Omega))\right)^n
\]

Let $\mathbf{E}(\Omega)\equiv E(\Omega)^n$ be the space
\[
\mathbf{E}(\Omega)=\lbrace \mathbf{u}\in \mathbf{L}^2(\Omega): \nabla\cdot \mathbf{u}\in L^2(\Omega)\rbrace
\]

This is a Hilbert space when equipped with the inner product
\[
\langle \mathbf{u},\mathbf{v}\rangle_{\mathbf{E}(\Omega)}=\langle \mathbf{u},\mathbf{v}\rangle_{\mathbf{L}^2(\Omega)}
+\langle \nabla\cdot\mathbf{u},\nabla\cdot\mathbf{v}\rangle_{L^2(\Omega)}.
\]
 
 The subspace of divergence free fields is defined by,

 \[
 \mathbf{E}_0(\Omega)=\lbrace \mathbf{h}_1\in\mathbf{E}(\Omega):\nabla\cdot\mathbf{h}_1=0\rbrace
 \]

\bigskip

Let us define the observation operator
\[
\mathcal{M}:L(\Omega)^3\rightarrow L(\Omega)^2,
\]
by
\[
\mathcal{M}(\mathbf{u})=(u_1,u_2)\equiv\mathbf{U}.
\]
It is readily seen that  in these $L^2$ spaces $\mathcal{M}^*\mathbf{U}=(U_1,U_2,0)$.

\bigskip

Let $\mathbf{U}^0\in L(\Omega)^2$ be the given field. Let us define
\[
J:\mathbf{E}(\Omega)\rightarrow \mathbb{R}
\]
by
\[
J(\mathbf{u})=\frac{1}{2}\Vert \mathcal{M}(\mathbf{u})-\mathbf{U}^0 \Vert^2_{\mathbf{L}^2(\Omega),S}
\]
where
\[
\langle \mathbf{U},\mathbf{V}\rangle_{\mathbf{L}^2(\Omega),S} =\int_\Omega (S\, \mathbf{U})\cdot \mathbf{V}
\]
and $S$ is a symmetric positive definite matrix.

\bigskip

We consider the problem: 
\begin{center}
minimize $J(\mathbf{u})\quad $ subject to $\quad\nabla\cdot \mathbf{u}=0$.
\end{center}

\subsection{Existence of minima}
 
 In what follows we shall use freely basic results from minimization of convex functionals on Hilbert spaces. See Zeidler \cite{zeidler}.
 
 \bigskip
   
 It is apparent that $J$ is Gateaux differentiable (G-differentiable). Let $\mathbf{h}_1,\mathbf{k}_1\in\mathbf{L}^2(\Omega)$, it is readily seen that the first variation of $J$ at $\mathbf{u}$ in the direction $\mathbf{h}_1$ is
 
 \begin{equation}
 \delta J(\mathbf{u};\mathbf{h}_1)=\langle\mathbf{h}_1,\mathcal{M}^*S(\mathcal{M}\mathbf{u}-\mathbf{U}^0)\rangle_{\mathbf{L}^2(\Omega)},
 \label{V1JL2}
 \end{equation} 
 whereas the second variation is given by
 
 \begin{equation}
 \delta^2J(\mathbf{u};\mathbf{h}_1,\mathbf{k}_1)=\langle \mathbf{h}_1,\mathcal{M}^*S\mathcal{M}\mathbf{k}_1\rangle_{\mathbf{L}^2(\Omega)}.
 \label{V2JL2}
 \end{equation}
   
\bigskip
   
We obtain the Taylor formula
\begin{equation}
 J(\mathbf{v})= J(\mathbf{u})+  \delta J(\mathbf{u};\mathbf{v}-\mathbf{u})
 +\frac{1}{2}\langle \mathbf{v}-\mathbf{u},\mathcal{M}^*S\mathcal{M}\mathbf{v}-\mathbf{u}\rangle_{\mathbf{L}^2(\Omega)}.
 \label{TalrJ}
 \end{equation}
 
 Note that for all $\mathbf{u},\mathbf{h}_1\in \mathbf{L}^2(\Omega)$
 
 \[
  \delta^2J(\mathbf{u};\mathbf{h}_1)\equiv  \delta^2J(\mathbf{u};\mathbf{h}_1,\mathbf{h}_1)\geq 0.
 \]
 
 Consequently $J$ is convex. Equivalently, for $\mathbf{v}\in\mathbf{L}^2(\Omega)$ we have
 
 \begin{equation*}
 J(\mathbf{v})\geq J(\mathbf{u})+ \delta J(u;\mathbf{v}-\mathbf{u}).
\end{equation*}

 \bigskip
 
 We are led to
 
 \bigskip
 
\noindent\textbf{Proposition. }If $N$ is a bounded, closed and convex subset of $\mathbf{L}^2(\Omega)$  the convex functional
\[
J:\mathbf{L}^2(\Omega)\, \to \mathbb{R}
\]
has a minimum.

\bigskip

We observe that $\mathbf{u}^0=(U^0_1,U^0_2,0)$ is a trivial minimum in $\mathbf{L}^2(\Omega)$ not necessarily in 
$\mathbf{E}_0(\Omega)$.  A correction is constructed below.
\section{The adjusted field by line search in $\mathbf{E}(\Omega)$}

For the computations that follow the next technical lemma for integration by parts is required.

\bigskip

\noindent\textbf{Lema 1. }Let $\Omega\subset\mathbb{R}^d,d\geq 2,$ be a bounded Lipschitz domain with boundary $\Gamma$. Then for all $\lambda\in H^1(\Omega)$ and $\mathbf{h}_1\in\mathbf{E}(\Omega)$,
\begin{equation}
\int_\Omega \lambda\nabla\cdot\mathbf{h}_1\, dx = \int_{\partial\Omega}\lambda\mathbf{h}_1\cdot\nu\, d\Gamma
-\int_\Omega\nabla\lambda\cdot\mathbf{h}_1\, dx.
\label{GenGreen}
\end{equation}
Where 
\[
\int_{\partial\Omega}\lambda\mathbf{h}_1\cdot\nu\, d\Gamma
\]
is well defined in the sense of the generalized trace with
\[
\mathbf{h}_1\cdot\nu\in H^{-1/2}(\Gamma),\quad\lambda\in H^{1/2}(\Gamma).
\]

\noindent\textbf{Proof. } (Lemma II.1.2.3 in  Sohr \cite{Sohr})

\bigskip


\subsection{The adjusted field}

Let $f(t)$ be the scalar quadratic function
\[
f(t)=J(\mathbf{u}_c+t\mathbf{p}).
\]
where  $\mathbf{u}_c$ is the base field and $\mathbf{p}$ a descendent direction. By (\ref{TalrJ}) we have
\begin{equation}
f(t)=J(\mathbf{u}_c)+t\langle \mathbf{p},\mathcal{M}^*S(\mathcal{M}\mathbf{u}_c-\mathbf{U}^0)\rangle_{\mathbf{L}^2(\Omega)}+
\frac{t^2}{2}
\langle \mathbf{p},\mathcal{M}^*S\mathcal{M}\mathbf{p}\rangle_{\mathbf{L}^2(\Omega)}
\label{fline}
\end{equation}

Consequently, for nontrivial $\mathbf{p}$ there exists a unique $t_c$ such that 
\[
\mathbf{u}_+=\mathbf{u}_c+t_c \mathbf{p}.
\]
is the minimizer along the line.

\bigskip

We call $\mathbf{u}_+$ the \emph{adjusted field} for the field $\mathbf{u}_c$ at direction $\mathbf{p}$.
\subsection{Incomplete data}

Let us consider $J$ restricted to $\mathbf{E}(\Omega)$. To complete the inner product in $\mathbf{E}(\Omega)$ we consider $\lambda\in H^2(\Omega)$ and write

\begin{equation}
\begin{array}{rcl}
f(t) & = & J(\mathbf{u}_c)+t\langle \mathbf{p},\mathcal{M}^*S(\mathcal{M}\mathbf{u}_c-\mathbf{U}^0)-\nabla\lambda\rangle_{\mathbf{E}(\Omega)}+
\frac{t^2}{2}
\langle \mathbf{p},\mathcal{M}^*S\mathcal{M}\mathbf{p}\rangle_{\mathbf{L}^2(\Omega)} \\
 & & \\
  & & + t(\langle \mathbf{p},\nabla\lambda\rangle_{\mathbf{L}^2(\Omega)}-
  \langle \nabla\cdot\mathbf{p},\nabla\cdot(\mathcal{M}^*S(\mathcal{M}\mathbf{u}_c-\mathbf{U}^0)-\nabla\lambda)\rangle_{\mathbf{L}^2(\Omega)})
 \end{array}  
\label{fline}
\end{equation}

For the multiplier we require
\begin{equation}
  \label{eq:poisson}
\Delta\lambda = \nabla\cdot(\mathcal{M}^*S(\mathcal{M}\mathbf{u}_c-\mathbf{U}^0)).
\end{equation}

We seek divergence free directions, $\nabla\cdot\mathbf{p}=0$, thus Lemma 1, implies
\[
\langle \mathbf{p},\nabla\lambda\rangle_{\mathbf{L}^2(\Omega)} = \int_\Gamma \lambda\mathbf{p}\cdot\nu d\Gamma.
\]

If $\lambda$ and $\mathbf{p}$ are chosen so that 
\begin{equation}
\label{integrall}
\int_\Gamma \lambda\mathbf{p}\cdot\nu d\Gamma = 0,
\end{equation}
then (\ref{fline}) becomes
\[
f(t) =J(\mathbf{u}_c)+t\langle \mathbf{p},\mathcal{M}^*S(\mathcal{M}\mathbf{u}_c-\mathbf{U}^0)-\nabla\lambda\rangle_{\mathbf{E}(\Omega)}+
\frac{t^2}{2} \langle \mathcal{M}\mathbf{p},S\mathcal{M}\mathbf{p}\rangle_{\mathbf{L}^2(\Omega)}.
\]

Hence, if the first two components of $\mathbf{p}$ form a nontrivial vector, there is a unique minimizer, which yields the 
adjusted field.

\bigskip

Let $\mathbf{p}$ be the steepest descent direction, namely,
\[
\mathbf{p}=-(\mathcal{M}^*S(\mathcal{M}\mathbf{u}_c-\mathbf{U}^0)-\nabla\lambda),
\]
and append one of the following boundary conditions to cancel the integral term (\ref{integrall}), in $f(t)$,

\begin{equation}
  \label{eq:bcond}
      (\mathcal{M}^*S(\mathcal{M}\mathbf{u}_c-\mathbf{U}^0)-\nabla\lambda)\cdot\nu=0  \textnormal{  or  }   \lambda = 0  \textnormal { for } \lambda \textnormal{ on }  \Gamma
\end{equation}

then 
\[
f(t) = J(\mathbf{u}_c)-t\langle \mathbf{p},\mathbf{p}\rangle_{\mathbf{L}^2(\Omega)}+
\frac{t^2}{2}
\langle S\mathcal{M}\mathbf{p},\mathcal{M}\mathbf{p}\rangle_{\mathbf{L}^2(\Omega)}.
\]

The adjusted field is given by
\[
u_+=u_c-t_c(\mathcal{M}^*S(\mathcal{M}\mathbf{u}_c-\mathbf{U}^0)-\nabla\lambda).
\]
Where
\[
t_c=\frac{\langle \mathbf{p},\mathbf{p}\rangle_{\mathbf{L}^2(\Omega)}}
{\langle S\mathcal{M}\mathbf{p},\mathcal{M}\mathbf{p}\rangle_{\mathbf{L}^2(\Omega)}}.
\]

\subsection{The Sasaki's adjusted field}

In the classical Sasaki's approach,  the functional
\[
J(\mathbf{u})=\int_\Omega \left[ \alpha_1^2(u_1-u_1^0)^2+\alpha_2^2(u_2-u_2^0)^2+\alpha_3^2(u_3-u_3^0)^2\right]dV
\]
is considered for correction. The initial vertical velocity $u_3^0$ is usually set to zero.

\bigskip

We address this problem in our general setting. The $L^2$ space is weighted by $S$, a $3\times 3$ symmetric positive definite matrix. 

Assume that a complete initial field $\mathbf{u}^0\in (L^2(\Omega))^3$ is given.
The corresponding observation operator is
\[
\mathcal{M}:\mathbf{L}(\Omega)^3\rightarrow \mathbf{L}(\Omega)^3,
\]
given by
\[
\mathcal{M}(\mathbf{u})=\mathbf{u}.
\]
The identity operator.

In this case we have
\[
f(t) = J(\mathbf{u}_c)+t\langle \mathbf{p},S(\mathbf{u}_c-\mathbf{u}^0)\rangle_{\mathbf{L}^2(\Omega)}+
\frac{t^2}{2}
\langle \mathbf{p},S\mathbf{p}\rangle_{\mathbf{L}^2(\Omega)} 
\]

Completing the inner product in $\mathbf{E}(\Omega)$ as before
\begin{equation}
f(t) = J(\mathbf{u}_c)+t\langle \mathbf{p},S(\mathbf{u}_c-\mathbf{u}^0)-\nabla\lambda\rangle_{\mathbf{L}^2(\Omega)}+
\frac{t^2}{2}\langle \mathbf{p},S\mathbf{p}\rangle_{\mathbf{L}^2(\Omega)}
+ t\langle \mathbf{p},\nabla\lambda\rangle_{\mathbf{L}^2(\Omega)}
\end{equation}

or

\begin{equation}
\begin{array}{rcl}
f(t) & = & J(\mathbf{u}_c)+t\langle S\mathbf{p},\mathbf{u}_c-\mathbf{u}^0-S^{-1}\nabla\lambda\rangle_{\mathbf{E}(\Omega)}+
\frac{t^2}{2}
\langle S\mathbf{p},\mathbf{p}\rangle_{\mathbf{L}^2(\Omega)} \\
 & & \\
  & & + t(\langle \mathbf{p},\nabla\lambda\rangle_{\mathbf{L}^2(\Omega)}-
  \langle \nabla\cdot (S\mathbf{p}),\nabla\cdot(\mathbf{u}_c-\mathbf{u}^0-S^{-1}\nabla\lambda)\rangle_{\mathbf{L}^2(\Omega)})
 \end{array}  
\end{equation}

For the multiplier $\lambda$ we now require
\begin{displaymath}
  \nabla\cdot(S^{-1}\nabla\lambda)= \nabla\cdot(\mathbf{u}_c-\mathbf{u}^0). 
\end{displaymath}

with boundary conditions 
\begin{displaymath}
  (\mathbf{u}_c-\mathbf{u}^0-S^{-1}\nabla\lambda)\cdot\nu=0 \textnormal{ or }  \lambda = 0,  \textnormal { for } \lambda \textnormal{ on }  \Gamma.
\end{displaymath}


Letting
\[
\mathbf{p}=-(\mathbf{u}_c-\mathbf{u}^0-S^{-1}\nabla\lambda),
\]
be the steepest descent direction, we obtain the adjusted field
\[
u_+=u_c-t_c(\mathbf{u}_c-\mathbf{u}^0-S^{-1}\nabla\lambda)
\]

It follows at once that 
\[
f(t) = J(\mathbf{u}_c)-t\langle S\mathbf{p},\mathbf{p}\rangle_{\mathbf{L}^2(\Omega)}+
\frac{t^2}{2}
\langle S\mathbf{p},\mathbf{p}\rangle_{\mathbf{L}^2(\Omega)}
\]
so that 
\[
t_c=1.
\]

This leads to
\begin{equation}
  \label{eq:sasaki}
  u_+=\mathbf{u}^0+S^{-1}\nabla\lambda.
\end{equation}
Hence, if $\mathbf{u}_c=0$ we obtain the Sasaki's solution.

\section{Boundary Conditions and RBF Discretization}

\begin{figure}[t]
  \centering
  \includegraphics[width=10.5cm,height=4cm]{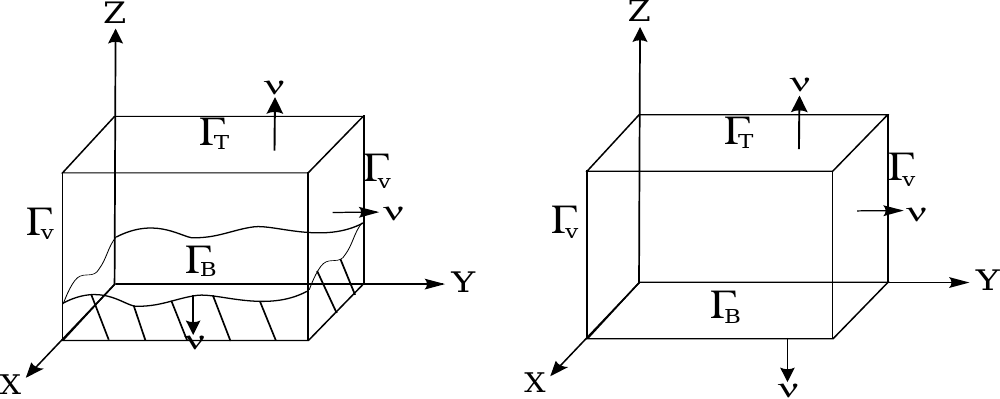}
 \caption{Bounded domain with and without topography.}
 \label{figu:bdomain} 
\end{figure}

For wind field recovery a bounded domain is considered as shown in Figure \ref{figu:bdomain}. An irregular bottom boundary models the topography of the terrain. With our formulation it is easily handled.

\subsection{Boundary Conditions}
\label{subsection:bc}

Given an initial field $\mathbf{U}^0$, a correction is through a line search in the direction of a Lagrange multiplier which satisfies the Poisson Boundary Value Problem (PBVP), given by (\ref{eq:poisson}) and (\ref{eq:bcond}). The latter clearly shows appropriate  boundary conditions, namely:

 \begin{itemize}
\item \textbf{\textit{Flow-Through}}.   The flux related to the vector field changes due to the mass balance through the boundary $\Gamma$. Thus $\nabla \lambda  \ne \mathbf{0}$, implying that one of the components (usually the normal component)  of the adjusted vector field is nonzero but unknown. To fulfill  (\ref{eq:bcond}) $\lambda = \mathbf{0}$ applies.


\item\textbf{\textit{No-Flow-Through}}.  This condition corresponds to the topography region $\Gamma_B$. There is no vertical component, hence $\textbf{u}_c \cdot  \nu = 0$.  

\item \textbf{ \textit {Dirichlet}}. In cases where the vector field is known on boundary $\Gamma$ (denoted by $\mathbf{u}_\Gamma$), we  set $\textbf{u}_c \cdot \nu  = \mathbf{u}_\Gamma\cdot \nu. $ For instance, in vertical borders $\Gamma_V$  the available information is the observed field $\mathbf{U}^0$. We require $\mathbf{u}_c\cdot\nu = \mathbf{U}^0\cdot \nu$.  

 \end{itemize}
 
\label{sec:na}

\subsection{RBF Asymmetric Collocation}

The asymmetric collocation method (\cite{kansa}), is the most  widely use technique for solving PDEs problems. The nodes are divided into interior and boundary nodes, the ansatz is build on each subsets of nodes and it is replaced by the analytic solution in the PDE and its boundary operators. This generates an algebraic system whose solution gives the coefficients of the approximated solution. 
Let $\Omega \subset \RR^n$ be an open bounded domain, let $X =\{ \bf\mathbf{c}_i \}_{i=1}^{N}$ a set of nodes  which are called centers. Define now the following PDE problem

\begin{align}
  \label{align:edpk}
  \begin{cases}
    \mathcal{L} \mathbf{u}(\mathbf{x}) = f(\mathbf{x}) & \mathbf{x} \in \Omega \\
    \mathcal{B} \mathbf{u}(\mathbf{x}) = g(\mathbf{x}) & \mathbf{x} \in \partial \Omega = \Gamma,
  \end{cases}
\end{align}Let the following radial bases function ansatz given by:
\begin{equation}\label{stencil}
  \hat{\mathbf{u}}(\mathbf{x}) :=  \sum_{j=1}^{N} \beta_j \phi(\|\mathbf{x}-\mathbf{c}_j\|_2)= \sum_{j=1}^{N} \beta_j \phi(r_j) 
\end{equation}
where  $\|\cdot\|_2 $ is the Euclidean norm. The collocation method is given by the following linear system   $G\mathbf{\beta} = b$, where $G$  named Gram's matrix is 
\begin{align}
  \label{eq:GM}
  (G)_{i,j} &=
  \mathcal{L}{\phi(\|\mathbf{x}-\mathbf{c}_j\|)}|_{\mathbf{x}= \mathbf{x}_i}  \quad \textnormal{ for } \mathbf{x}_i\in \Omega    \nonumber \\   
  (G)_{i,j} &=
  \mathcal{B}{\phi(\|\mathbf{x}-\mathbf{c}_j\|)}|_{\mathbf{x}=\mathbf{x}_i}   \quad \textnormal{ for } \mathbf{x}_i  \in  \Gamma \nonumber
\end{align} and
\begin{displaymath}
  b=\left[
  \begin{matrix}
     b_I \\
     b_{\Gamma} \\
   \end{matrix}  
 \right] 
\end{displaymath}
with
\begin{align*}
  b_I &= [f(\mathbf{x}_i)]^T \textnormal{ for } \mathbf{x}_i\in \Omega \\
  b_\Gamma &= [f(\mathbf{x}_i)]^T \textnormal{ for } \mathbf{x}_i\in \Gamma
\end{align*}

\section{Numerical examples}
\label{sec:numerical}

In the following examples we illustrate the application of the line search method above. We analyse  its performance when recovering the unknown vertical component of zero divergence fields. In all examples we use Dirichlet data. 

For the numerical approximation we shall use the ansatz (\ref{stencil}), using inverse multiquadrics,
\begin{displaymath}
  \phi(r) = \frac{1}{\sqrt{(1+(rc)^2)}}
\end{displaymath}
where  $c$ is the shape parameter, whose value is related to the spectral  convergence of the approximation (\cite{fornberg}).

\begin{exmp}
  \label{example2}

\begin{figure}[t] 

    \centering
    \includegraphics[scale=.7]{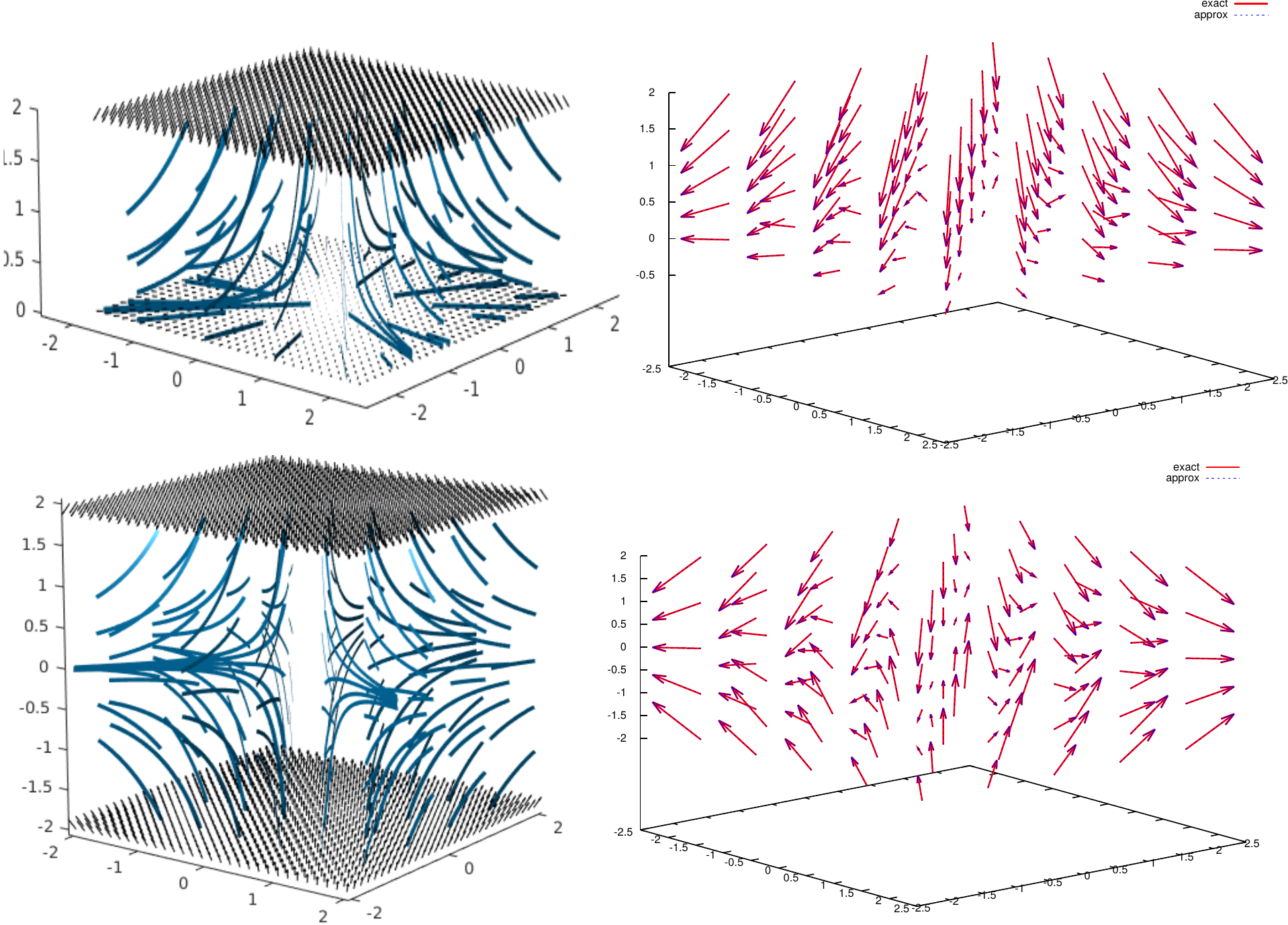}
    \caption{{\bf Example 5.1. }(Left)Vector field   $\mathbf{u}(x,y,z) = (x,y,-2z)$ on $(-2,2)\times(-2,2)\times(0,2)$ and $(-2,2)\times (-2,2)\times (-2,2)$. (Right) Vector field approximation taking $\mathbf{u}^0=(x,y,0)$, the arrows {\color{red} $-\!\!\!\!\blacktriangleright$}  and   {\color{blue} $\cdots\!\!\!>$}, represents respectively exact and approximated vector field.}
    \label{figure:fig3y4} 
\end{figure}
Let us consider an extension of the 2D field in Cervantes et al \cite{cervantesetal}, namely
 $\mathbf{u}(x,y,z)=(x,y,-2z)$ on the  domains $\Omega_1=(-2,2)\times(-2,2)\times(0,2)$ and $\Omega_2=(-2,2)\times(-2,2)\times(-2,2)$ (which are shown on the left side of figure  \ref{figure:fig3y4} in a stream ribbon style).  Assume that  $\mathbf{u}^0(x,y,z)=(x,y,0)$, so that, in   $\Omega_1$ (domain with topography), we impose a No-Flow-Thrown boundary condition on $\Gamma_B$, while in $\Omega_2$ a Dirichlet boundary condition is used on $\Gamma_B$ (see figure \ref{figu:bdomain}).

The results are shown in table \ref{table:tabla2} and in the right hand side of figure \ref{figure:fig3y4}. In the table \ref{table:tabla2}, only the results for $\Omega_1$ are displayed.  The corresponding results  for $\Omega_2$ are nearly equivalent.  The relative error  and the approximated divergence
, decreases as the number of centers increases. As expected,  the  condition number of  Gram matrix increases. This is the well know and  classic behavior of the schemes based on radial basis functions, the so called \textit{Schaback Uncertainty Principle} \cite{schaback}.

 \begin{table}[h]
\centering
\begin{tabular}{c c c c c}
\hline
 $N$ & $c$  &  $\kappa(G) $& $\nabla \cdot u_+$ & $\|u_{+}-\mathbf{u}\|_{2}/\|\mathbf{u}\|_2 $   \\ 
\hline 
  27   & 0.001 & 2.252627e+18   & 4.679924e-05 &  2.707875e-05    \\ 
  125 & 0.001 & 5.801561e+19 & 1.088109e-06 & 5.342928e-06   \\ 
  512 & 0.001 & 1.195701e+20 & 7.873625e-08  &  6.961732e-08   \\ 
\hline \\
\end{tabular}
    \caption{{\bf Example 5.1.} Vector field approximation ($\mathbf{u}_+$)  of   $\mathbf{u}(x,y) = (x,y,-2z)$ on $(-2,2)\times(-2,2)\times (0,2)$ using a set of equidistant centers. }
\label{table:tabla2}
\end{table}
  \end{exmp}

\begin{exmp}
\label{exmp:helicoidal_topo}
\begin{figure}[t] 
    \label{figure:bdomain} 
    \centering
    \includegraphics[scale=.75]{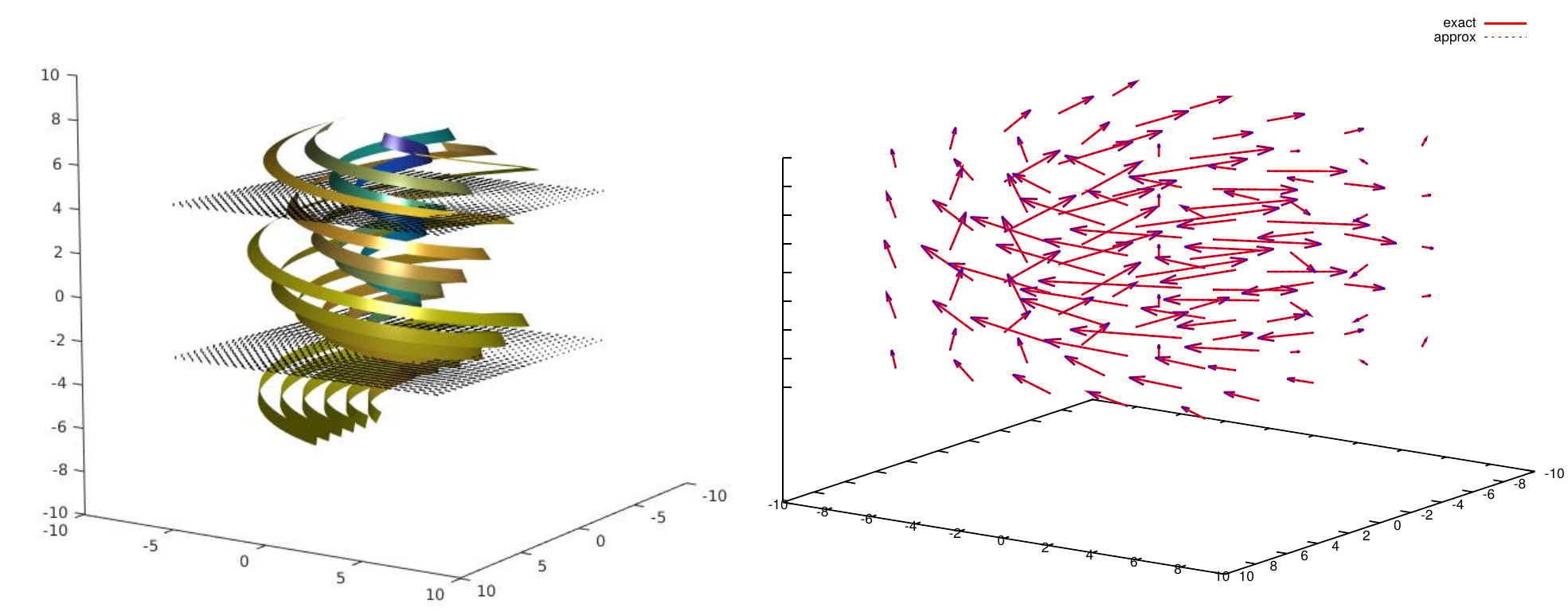}
    \caption{ {\bf Example 5.2.} (Left)  Vector field  $\mathbf{u}(x,y,z) = (2ye^{\frac{-(x^2+y^2+z^2)}{49}},-2xe^{\frac{-(x^2+y^2+z^2)}{49}},1)$ on $(-7,7)\times(-7,7)\times(-7,7)$.  (Right) Vector field approximation taking taking  $\mathbf{u}^0(x,y,z) = (2ye^{-\frac{(x^2+y^2+z^2)}{49}},-2xe^{- \frac{(x^2+y^2+z^2)}{49}},0 )$.  The arrows {\color{red} $-\!\!\!\!\blacktriangleright$}  and   {\color{blue} $\cdots\!\!\!>$}, represents respectively exact and approximated vector field.}
\label{figure:fig7y8}
\end{figure}
In this example, we approximate the (nontrivial) vortex type vector field in Benbourhim et al \cite{najib}. It is shown on the left side of figure  \ref{figure:fig7y8} in a stream ribbon style, and  is  given by
  \begin{displaymath}
   \mathbf{u}(x,y,z) = (2ye^{-\frac{(x^2+y^2+z^2)}{49}} ,-2xe^{- \frac{(x^2+y^2+z^2)}{49}},1 )
 \end{displaymath}
  with  $\Omega = (-7,7)\times (-7,7)\times(-7,7)$. The initial field is taken as
\begin{displaymath}
    \mathbf{U}^0(x,y,z) = (2ye^{-\frac{(x^2+y^2+z^2)}{49}} ,-2xe^{- \frac{(x^2+y^2+z^2)}{49}},0 )
  \end{displaymath}
In this case, we do not consider the topography and assume that vector field in the upper and lower regions are known. That is, Dirichlet boundary conditions for $\Gamma_T$ and $\Gamma_B$ are considered. The results are shown in table \ref{table:table3} and in the right hand side of figure \ref{figure:fig7y8}.

\begin{table}[h]
\centering
\begin{tabular}{c  c c c c}
\hline 
 $N$ & $c$  &   $\kappa(G) $& $\nabla \cdot u_+$ & $\|u_{+}-\mathbf{u}\|_{2}/\|\mathbf{u}\|_2 $    \\ 
\hline 
   27  & 0.01 &   6.468932e+09   &  -5.553522e-06  &  4.879887e-03    \\ 
  125 & 0.01  & 5.736571e+18 &  -6.975779e-03 &  2.968737e-05   \\ 
  512 & 0.01  & 1.057278e+20 & -5.411860e-03  &  1.226800e-07    \\ 
  \hline \\
\end{tabular}

    \caption{{\bf Example 5.2.}  Vector field approximation ($\mathbf{u}_+$)  of   $ (2ye^{-\frac{(x^2+y^2+z^2)}{49}} ,-2xe^{- \frac{(x^2+y^2+z^2)}{49}},1 )$ on $(-7,7)\times(-7,7)\times (-7,7)$ using a set of equidistant centers. }
\label{table:table3}
\end{table}

 \end{exmp}

\begin{exmp}
\label{exmp:helicoidal_topo}
\begin{figure}[t] 
    \label{figure:bdomain} 
    \centering
    \includegraphics[scale=.75]{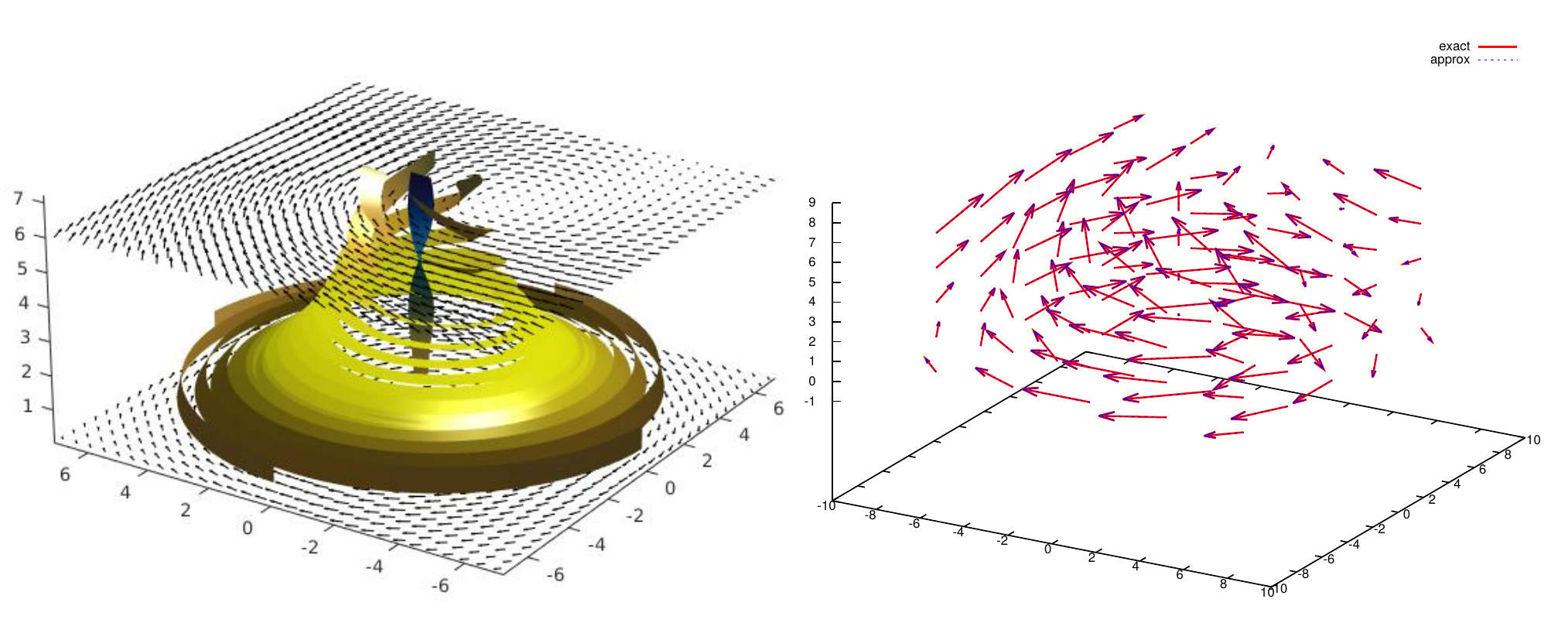}
    \caption{ {\bf Example 5.3.} (Left) Vector field  $\mathbf{u}(x,y,z) = (2ye^{\frac{-(x^2+y^2+z^2)}{49}}-\varepsilon \frac{xz}{2},-2xe^{\frac{-(x^2+y^2+z^2)}{49}}-\varepsilon \frac{xz}{2},\varepsilon \frac{z^2}{2})$ on $(-7,7)\times(-7,7)\times(0,7)$, for  $\varepsilon = 0.1$. (Right) Vector field approximation,  taking  $\mathbf{u}^0(x,y,z) = (2ye^{-\frac{(x^2+y^2+z^2)}{49}} - \varepsilon \frac{ x z}{2},-2xe^{- \frac{(x^2+y^2+z^2)}{49}}-\varepsilon \frac{yz}{2},0 )$. The arrows {\color{red} $-\!\!\!\!\blacktriangleright$}  and   {\color{blue} $\cdots\!\!\!>$}, represents respectively exact and approximated vector field.}
 \label{figure:fig5y6}
\end{figure}
In this last example we consider a variant of Example 5.2 to have a topography in the $xy$ plane. Consequently the No-Flow-Through boundary condition is imposed on $\Gamma_B$. 

The vector field is given by, 

\begin{displaymath}
   \mathbf{u}(x,y,z) = (2ye^{-\frac{(x^2+y^2+z^2)}{49}} - \varepsilon \frac{ x z}{2},-2xe^{- \frac{(x^2+y^2+z^2)}{49}}-\varepsilon \frac{yz}{2},\varepsilon \frac{z^2}{2} )
 \end{displaymath}
  with  $\Omega = (-7,7)\times (-7,7)\times(0,7)$ and $\varepsilon  > 0 \in \RR^+$. As before, the initial field is, 
  \begin{displaymath}
    \mathbf{u}^0(x,y,z) = (2ye^{-\frac{(x^2+y^2+z^2)}{49}} - \varepsilon \frac{ x z}{2},-2xe^{- \frac{(x^2+y^2+z^2)}{49}}-\varepsilon \frac{yz}{2},0 ).
  \end{displaymath}
In the left hand side of figure \ref{figure:fig5y6},  we  show this vector filed for $\varepsilon = 0.1$  in a stream ribbon style. The corresponding results are shown in table \ref{table:table4} and in the right hand side of figure \ref{figure:fig5y6}.
\begin{table}[h]
\centering
\begin{tabular}{c c c c c c}
\hline 
 $N$ & $c$  & $\epsilon$ &  $\kappa(G) $& $\nabla \cdot u_+$ & $\|u_{+}-u\|_{2}/\|u\|_2 $    \\ 
\hline  
\hline
   27  & 0.01 &0.1 &  1.508177e+13   & 2.633422e-03 &  1.463895e-01    \\ 
  125 & 0.01 &0.1 & 2.632883e+21  & 4.063153e-03 & 4.732374e-04   \\ 
  512 & 0.01 &0.1 & 8.866826e+21 &1.365358e-02  & 5.872183e-05    \\ 
\hline \\
\end{tabular}
    \caption{{\bf Example 5.3.}  Vector field approximation ($\mathbf{u}_+$)  of   $\mathbf{u}(x,y) =  (2ye^{-\frac{(x^2+y^2+z^2)}{49}} - \varepsilon \frac{ x z}{2},-2xe^{- \frac{(x^2+y^2+z^2)}{49}}-\varepsilon \frac{yz}{2},\varepsilon \frac{z^2}{2} )$ on $(-7,7)\times(-7,7)\times (0,7)$ using a set of equidistant centers. }
\label{table:table4}
\end{table}

 \end{exmp}

\section{Conclusion} 

In this article we have introduced a mass consistent variational approach to approximate 3D vector fields from a set of prescribed scattered data values, whose vertical component is unknown. This problem is of great interest for meteorological applications and has been treated in the literature by several authors; see Ratto et al \cite{ratto}, for a recent review.

More precisely a first contribution of this article is that our approach, based on the adjoint method, let us build the adjusted field, as the unique solution along an appropriately chosen descent direction. A full analysis of the existence and uniqueness  of the solution is proved within the proper Hilbertian setting. Boundary conditions for the adjoint equations arises naturally from the variational formulation in such a way that it is possible to build them according to the physical properties of the field. This point generalizes Sasaki's pioneer approach, \cite{sasaki}, which imposes in an heuristic way, particular boundary conditions. We prove in fact that Sasaki's method is a particular case of our approach.

Several works, which uses classic techniques, like finite elements or finite differences, have appeared in the literature, to solve this problem, (see, \cite{ratto}).

In order to avoid the computational expansive 3D mesh generation, some Radial Basis Functions (RBF) methods have been formulated;  See Pepper, Rasmussen \& Fyda \cite{Pepperetal}, and Cervantes et al \cite{cervantesetal}.

We recall the work of Benbourhim and Bouhamidi \cite{najib}, that developed a smoothing algorithm for a given set of physically constrained prescribed data of vector fields in two and three dimensions. Although this is an important method, no boundary conditions can be imposed by this approach in an open bounded set in tow or three dimensions.

In our work a radial mesh free asymmetric collocation technique is used to discretize the adjoint equations in an open bounded domain. Inverse multiquadric kernels are used, thus providing exponential rate of convergence. Numerical experiments were designed to prove different boundary conditions, namely, when the field is tangential to the terrain and when no terrain is considered. 
In all cases, the numerical errors and the approximated value of the divergence are excellent for a small number of data. 

Several points remains unsolved within the scope of  this work.  Among them, we mention its possible generalization to approximate vector fields having different physical properties, like zero rotational, or elasticity constrains. Also, in order to approximate real life data, this method could be restated as a smoothing approximant. Further work is in progress in these directions.

\section{Acknowledgements} 

The authors would like to acknowledge to ECOS-NORD project number \newline 000000000263116/M15M01 for financial support during this research. C. Gout thanks the M2NUM project which is co-financed by the European Union with the European regional development fund (ERDF, HN0002137) and by the Normandie Regional Council."

\end{document}